\documentclass{article}
\def\dil{\Lam}

\def\capset{\Gamma}
\def\dualset{\capset^\ast}

\newtheorem{theorem}{Theorem}
\def\beaN{\setlength{\arraycolsep}{0.0em}\begin{eqnarray*}}
\def\eeaN{\end{eqnarray*}\setlength{\arraycolsep}{5pt}}
\def\bea{\setlength{\arraycolsep}{0.0em}\begin{eqnarray}}
\def\eea{\end{eqnarray}\setlength{\arraycolsep}{5pt}}

\def\be{\begin{equation}}
\def\ee{\end{equation}}

\def\dm{n}

\def\Zd{\ZZ^\dm}

\def\RR{\mathbb{R}}
\def\ZZ{\mathbb{Z}}

\def\bks{\backslash}

               \def\Lam{\Lambda}

\def\ome{\omega}                

\usepackage{cite}
\usepackage{graphicx}
\usepackage{amssymb}
\usepackage{enumerate}
\usepackage{bm,bbm}
\usepackage{spconf,amsmath}

\title{COMMITTEE ALGORITHM: AN EASY WAY TO CONSTRUCT WAVELET FILTER BANKS}
\name{Youngmi Hur
}
\address{Department of Applied Mathematics and Statistics\\
The Johns Hopkins University, Baltimore, MD 21218}
\begin{document}
%\ninept
%
\maketitle
\begin{abstract}
Given a lowpass filter, finding a dual lowpass filter is an essential step in constructing non-redundant wavelet filter banks. Obtaining dual lowpass filters is not an easy task. In this paper, we introduce a new method called committee algorithm that builds a dual filter straightforwardly from two easily-constructible lowpass filters. It allows to design a wide range of new wavelet filter banks. An example based on the family of Burt-Adelson's 1-D Laplacian filters is given.\end{abstract}
\begin{keywords}
Laplacian pyramid, non-redundant filter bank, polyphase representation, wavelet filter bank
\end{keywords}
\section{Introduction}
\label{sec:intro}
A filter bank (FB) consists of the analysis bank and the synthesis bank, which are collections of, say $p$, filters linked by downsampling and upsampling operators, respectively \cite{SN}. 
A FB is typically referred to as a {\it wavelet} FB if each of its analysis and synthesis banks has exactly one lowpass filter and the rest of them are all highpass filters. Designing non-redundant wavelet FBs is an important problem since it leads to the construction of biorthogonal (or Riesz) wavelet bases under well-understood constraints \cite{Ma1,Me,SN}.

Such a design problem is often reduced to the problem of extending a matrix with Laurent polynomial entries \cite{LLS}. In this approach, an essential step is to find a dual lowpass filter $d$ of a given lowpass filter $h$ so that (i) $h$ and $d$ are biorthogonal, (ii) $d$ has positive accuracy\footnote{In designing non-redundant {\it wavelet} FBs, having lowpass filters with {\it positive} accuracy is essential as the positive accuracy guarantees that all the other filters are highpass \cite{CHR}.} (cf. see Section 2 for definition), and optionally (iii) $d$ has some other properties desirable for the specific design. Although a filter that is only biorthogonal to $h$, or a filter that has only positive accuracy can be obtained readily, finding the filter that satisfies both properties (i) and (ii) at the same time is not easy.
The filters with biorthogonality property but not with positive accuracy lead to FBs but not to wavelet FBs.

In this paper we show that the dual filter $d$ can be obtained by first finding two lowpass filters $f$ and $g$ that have a ``partial responsibility'': (a) $f$ is biorthogonal to $h$ but it needs not have positive accuracy and (b) $g$ has positive accuracy but it needs not be biorthogonal to $h$. We show that given such a pair of lowpass filters $f$ and $g$, the dual filter $d$ that satisfies the above conditions (i) and (ii) simultaneously can be obtained straightforwardly. Since finding the pair $(f,g)$ that satisfies (a) and (b) independently is much easier than finding the dual $d$ directly, our method provides an easier algorithm to construct non-redundant wavelet FBs. We refer to this algorithm as the {\it committee algorithm}. 

In order to develop our committee algorithm we use the Laplacian pyramid (LP) algorithm and its polyphase representation. We briefly review these along with other relevant concepts in Section 2. Our main results are presented in Section 3 together with an example illustrating our findings. We summarize our results in Section 4.

\section{LP and its polyphase representation}

We consider only the FBs with finite impulse response (FIR) filters that have the perfect reconstruction (PR) property. FIR FBs are useful as they provide fast algorithms.
We use $\dil$ to denote the $\dm\times\dm$ dilation matrix for sampling, where $\dm$ is the spatial dimension. The PR property holds only if $p\ge q:=|\det \dil|$, where $p$ is the number of filters in each bank as before. If $p=q$, the FB is called {\it non-redundant}; otherwise, it is called redundant. 

The Laplacian pyramid is introduced by Burt and Adelson \cite{BA}. The LP algorithm has analysis and synthesis processes. One begins the LP analysis process with two FIR lowpass filters: the compression and prediction filters. Using the compression filter $h$, one obtains the coarse coefficients
that approximate the input signal. The prediction
  filter $g$ is used for predicting the original signal from the coarse
coefficients. One computes the detail coefficients by subtracting the predicted signal from the input signal, and then stores them along with the coarse coefficients in place of the original signal.
The standard LP synthesis process recovers the original signal by reversing the above subtraction step. It is well-known that the LP representation is {\it redundant}; there are more coefficients after
the LP analysis process than those in the input signal. 

The polyphase decomposition in \cite{V} is a widely used
method of transforming a filter (or signal) into $q=|\det\dil|$
filters (or signals) running at the sampling rate $1/q$. Below we briefly review the polyphase representation of
LP and refer \cite{DV2,Hur} for details. We let $\capset$ be a
complete set of representatives of the distinct cosets of the quotient
group $\Zd/\dil\Zd$ containing $0$, and $\dualset$ be a complete
set of representatives of the distinct cosets of
$2\pi(((\dil^T)^{-1}\Zd)/\Zd)$ containing $0$. Here $\dil^T$ is used to denote the transpose of $\dil$. Then both the sets $\capset$
and $\dualset$ have $q=|\det \dil|$ elements. We use 
$$\nu_0=0, \nu_1,\cdots, \nu_{q-1}$$
to denote the elements of $\capset$.

The (polyphase) LP filters are 
\beaN
{\tt G}(z)&{\,:=\,}&[G_{\nu_0}(z),G_{\nu_1}(z),\dots,G_{\nu_{q-1}}(z)]^T,\\
{\tt H}(z)&{:=}&[H_{\nu_0}(z),H_{\nu_1}(z),\dots,H_{\nu_{q-1}}(z)],
\eeaN
where $G_\nu$ and $H_\nu$ are the $z$-transforms \cite{RZ} of $g_\nu$ and $h_\nu$, respectively, with
$$g_\nu(m):=g(\dil m+\nu), \quad h_\nu(m):=h(\dil m-\nu), \quad\forall m\in\Zd,$$
and  
$$
{\tt A}_{LP}(z)
:=\left[\begin{array}{c}
                {\tt H}(z)\\
                {\tt I}_q-{\tt
  G}(z){\tt H}(z)
               \end{array}\right],
%=:\left[\begin{array}{c}
%                {\tt H}(z)\\
%                \tilde{\tt A}_{LP}(z)
%               \end{array}\right].
$$
\begin{equation}
\label{eq:S0}
{\tt S}_0(z):=\left[\begin{array}{cc}
                {\tt G}(z) & {\tt I}_q
               \end{array}\right],
%=:\left[\begin{array}{cc}
%                {\tt G}(z) & \tilde{\tt S}_0(z)
%               \end{array}\right],
\end{equation}
are the (polyphase) analysis and trivial synthesis operators of LP, respectively. Here and below we denote the $m\times m$ identity matrix using ${\tt I}_m$.
Since $h$ and $g$ are FIR filters, the entries of the above operators are all Laurent polynomials in $z$.

In this paper we normalize the lowpass filter so that its coefficients sum up to $\sqrt{q}$. 
The LP compression filter is called {\it interpolatory} if $H_{0}(z)=1/\sqrt{q}$. 
We also assume that both LP filters have {\it positive accuracy}.
We recall that for an integer $m\ge 0$, the filter $h$ has {\it accuracy $m$} if the number of zeros of $H(e^{i\omega})$, the Fourier transform of $h$, at
$\omega\in\dualset\bks\{0\}$ is $m$ \cite{SN}. 

The pair $({\tt A}_{LP}(z), {\tt S}_0(z))$ of the above LP operators can be considered as the polyphase representation of a FB. 
The fact that LP is redundant can be seen also from these polyphase representations. The analysis operator ${\tt A}_{LP}(z)$ and the trivial synthesis operator ${\tt S}_0(z)$ clearly satisfy the perfect reconstruction property ${\tt S}_0(z){\tt A}_{LP}(z)={\tt I}_q$, but it is oversampled since ${\tt A}_{LP}(z)$
and ${\tt S}_0(z)$ are not square matrices. They are of size $(q+1)\times q$ and $q\times (q+1)$, respectively. 

For a given analysis operator ${\tt A}_{LP}(z)$, there exist infinitely many different synthesis operators that satisfy the perfect reconstruction property. It is well-known (\cite{BHF,DV2,Hur,LGT}) that the most general LP synthesis operator ${\tt S}_{LP}(z)$ is given as the form
\begin{equation}
\label{eq:generalS}
{\tt S}_{LP}(z):=\left[\begin{array}{cc}
                {\tt G}(z)+{\tt V}(z)B(z) & {\tt
                I}_q-{\tt V}(z){\tt H}(z)
               \end{array}\right],
\end{equation}
with $B(z):=1-{\tt H}(z){\tt G}(z)$, ${\tt V}(z):={\tt U}_c(z)+{\tt U}_d(z){\tt G}(z)$, where ${\tt U}_c(z)$ and ${\tt U}_d(z)$ are the
parameter matrices of size
$q\times 1$ and $q\times q$, respectively, consisting of the
Laurent polynomial entries. 
We note that the LP filters are biorthogonal if and only if $B(z)=0$. 

The synthesis operator ${\tt S}_0(z)$ in
(\ref{eq:S0}) corresponds to the case where ${\tt V}(z)=0$.
The LP pair $({\tt A}_{LP}(z), {\tt S}_0(z))$, with the trivial synthesis operator, is not a wavelet FB since none of the filters in the synthesis bank ${\tt S}_0(z)$ is highpass.

There are a few existing methods for designing wavelet FBs from LPs.
The wavelet frame theory was used in \cite{DV2,HR1} to 
reinterpret the LP process as a wavelet process, and the associated wavelet FB was constructed along the way. However the resulting wavelet FB is redundant and the filters are not necessarily FIR unless the LP filters are biorthogonal. 
The paper \cite{Hur} concerns about transforming the given LP to a critical representation (CLP). 
A CLP of a given LP is said to exist if there exists a $q\times (q+1)$ matrix ${\tt R}(z)$, called the reduction operator, whose entries are the Laurent polynomials in $z$ such that (i) ${\tt R}(z){\tt A}_{LP}(z)$ is invertible, and (ii) the first row of ${\tt R}(z){\tt A}_{LP}(z)$ is the same as ${\tt H}(z)$. We also recall that a CLP is called ECLP if the reduction operator takes a simple special form. The ECLP enables us to obtain a critical representation effortlessly by eliminating a redundant portion of the LP. In the same article, it is shown that a special class of ECLP called interpolatory is available whenever the LP compression filter satisfies 
\begin{equation}
\label{eq:ECLPthm}
H_{\nu_{\beta-1}}(z)=1/\sqrt{q},\quad \hbox{for some $\beta\in\{1,\dots,q\}$,}
%={1\over \sqrt{q}},\quad\hbox{ for some $\nu\in\capset$.}
\end{equation}
which includes the LP compression filter $h$ being interpolatory as a special case, is satisfied. It turns out that the interpolatory ECLP also provides a way to design non-redundant FIR wavelet FBs from LPs with the condition (\ref{eq:ECLPthm}).

\section{Committee Algorithm}

\subsection{Theory}

It is easy to see that a CLP exists if and only if there exists a non-redundant FIR FB $({\tt A}_{CLP}(z), {\tt S}_{CLP}(z))$ such that the first row of ${\tt A}_{CLP}(z)$ is the same as the given LP compression filter ${\tt H}(z)$. It is also easy to see that if such a non-redundant FIR FB exists, then there exists Laurent
polynomials $f_{\nu_0}(z), f_{\nu_1}(z),\cdots,f_{\nu_{q-1}}(z)$ such that 
\begin{equation}
\label{eq:coprime}
f_{\nu_0}(z)H_{\nu_0}(z)+f_{\nu_1}(z)H_{\nu_1}(z)+\cdots+f_{\nu_{q-1}}(z)H_{\nu_{q-1}}(z)=1.
\end{equation}
Since ${\tt F}(z):=[f_{\nu_0}(z),f_{\nu_1}(z),\cdots,f_{\nu_{q-1}}(z)]^T$ can be considered as the polyphase representation of a filter, say $f$, the above condition between ${\tt F}(z)$ and ${\tt H}(z)$ is equivalent to the biorthgonality of the filters $f$ and $h$. Furthermore it is easy to see that the filter $f$ whose polyphase representation satisfies (\ref{eq:coprime}) is necessarily lowpass (see, for example, \cite{Hur}). 

The next theorem shows that the above condition is also sufficient to obtain a wavelet FB and that there is a simple algorithm which produces the wavelet FB. 

\begin{theorem}
Suppose that the LP compression filter satisfies (\ref{eq:coprime}) for some Laurent
polynomials $f_{\nu_0}(z), f_{\nu_1}(z),\cdots,f_{\nu_{q-1}}(z)$. Then there exists a CLP with a wavelet FB.
\end{theorem}

{\sl Proof: }
{\sl 
Let ${\tt V}(z)=[f_{\nu_0}(z), f_{\nu_1}(z),\cdots,f_{\nu_{q-1}}(z)]^T$. Then, from $\det({\tt I}_{q}-{\tt V}(z){\tt H}(z))=1-{\tt H}(z){\tt V}(z)=0$ and ${\tt S}_{LP}(z){\tt A}_{LP}(z)={\tt I}_q$, and from the form of ${\tt S}_{LP}(z)$ given in (\ref{eq:generalS}), the rank of ${\tt I}_{q}-{\tt V}(z){\tt H}(z)$ has to be $q-1$. Thus, there are exactly $q-1$ linearly independent columns of ${\tt I}_{q}-{\tt V}(z){\tt H}(z)$ and the other column can be written as a linear combination of these $q-1$ columns. Suppose that the $k$-th column can be written as a linear combination of the rest of the columns. Then there exists a coefficient vector 
$${\tt c'}:=[c_1,\cdots,c_{k-1}]^T,\quad {\tt c''}:=[c_{k+1},\cdots,c_q]^T$$
such that $({\tt I}_{q}-{\tt V}(z){\tt H}(z))[{\tt c}',1,{\tt c}'']^T=0$. Thus the $(k+1)$-th column of 
$${\tt S}(z)={\tt S}_{LP}(z)\left[\begin{array}{cccc}
                1 & 0             & 0 & 0\\
                0 & {\tt I}_{k-1} & {\tt c'} & 0\\
                0 & 0 & 1 & 0\\
                0 & 0 & {\tt c''} & {\tt I}_{q-k}
               \end{array}\right]$$
is zero vector. Since ${\tt S}_{LP}(z){\tt A}_{LP}(z)={\tt I}_q$ and 
$$\left[\begin{array}{ccc}
               {\tt I}_{k-1} & {\tt c'} & 0\\
                0 & 1 & 0\\
                0 & {\tt c''} & {\tt I}_{q-k}
               \end{array}\right]\left[\begin{array}{ccc}
                {\tt I}_{k-1} & -{\tt c'} & 0\\
                0 & 1 & 0\\
                0 & -{\tt c''} & {\tt I}_{q-k}
               \end{array}\right]={\tt I}_{q},$$
if we define ${\tt S}_{CLP}(z)$ as the submatrix of ${\tt S}(z)$ obtained by deleting its $(k+1)$-th column, and ${\tt A}_{CLP}(z)$ as the submatrix of
$${\tt A}(z)=\left[\begin{array}{cccc}
                1 & 0             & 0 & 0\\
                0 & {\tt I}_{k-1} & -{\tt c'} & 0\\
                0 & 0 & 1 & 0\\
                0 & 0 & -{\tt c''} & {\tt I}_{q-k}
               \end{array}\right]{\tt A}_{LP}(z)$$
by deleting its $(k+1)$-th row, we have ${\tt S}_{CLP}(z){\tt A}_{CLP}(z)={\tt I}_q$. It is easy to see that the first row of ${\tt A}_{CLP}(z)$ is ${\tt H}(z)$. Hence the CLP with the FB $({\tt A}_{CLP}(z), {\tt S}_{CLP}(z))$ exists.
%Note that the reduction operator of this CLP is 
%$$
%\left[\begin{array}{cccc}
%                1 & 0             & 0 & 0\\
%                0 & {\tt I}_{k-1} & -{\tt c'} & 0\\
%                0 & 0 & -{\tt c''} & {\tt I}_{q-k}
%               \end{array}\right]
%$$

Let us now show that the FB $({\tt A}_{CLP}(z), {\tt S}_{CLP}(z))$ is a wavelet FB.
Let $d$ be the filter associated with the polyphase representation ${\tt G}(z)+{\tt V}(z)B(z)$. Since the FB is non-redundant, it suffices to show that the filters $h$ and $d$ have positive accuracy (see for example \cite{CHR}). Since $h$ is assumed to have positive accuracy, all we need to show is that $d$ has positive accuracy, which is equivalent to (cf.~Result 2 in \cite{Hur}) ${\tt G}({\tt 1})+{\tt V}({\tt 1})B({\tt 1})={1\over\sqrt{q}}{\tt 1}$ where ${\tt 1}=[1,\cdots,1]^T\in\RR^q$.
This is true since, no matter what the value ${\tt V}({\tt 1})$ is, we have 
\begin{eqnarray*}
&{ }&{\tt G}({\tt 1})+{\tt V}({\tt 1})B({\tt 1})
 = {\tt G}({\tt 1})+{\tt V}({\tt 1})(1-{\tt H}({\tt 1}){\tt G}({\tt 1}))\\
&=&{1\over\sqrt{q}}{\tt 1}+{\tt V}({\tt 1})(1-{1\over\sqrt{q}}{\tt 1}^T{1\over\sqrt{q}}{\tt 1})
={1\over\sqrt{q}}{\tt 1}
\end{eqnarray*}
where the fact that the LP filters $h$ and $g$ have positive accuracy is used for the second equality.}

\medskip
It is easy to see that if the LP compression filter satisfies the condition (\ref{eq:ECLPthm}) for the interpolatory ECLP, then the set of Laurent polynomials $\{f_{\nu}(z):\nu\in\capset\}$ satisfying (\ref{eq:coprime}) always exists as one can choose $f_{\nu_{\beta-1}}(z)=\sqrt{q}$ and $f_{\nu_{\alpha-1}}(z)=0$ for all $\alpha\ne \beta$. Hence our result here can be considered as a generalization of the interpolatory ECLP method in \cite{Hur}.

\subsection{Algorithm}

Theorem 1 provides a new method to design (non-redundant FIR) wavelet FBs.
In particular we obtain an algorithm to construct a dual lowpass filter from a given lowpass filter. 

\medskip
{\noindent\it Algorithm: Constructing a dual filter of a given lowpass filter}

{\noindent\bf Input:} $h$, a lowpass filter with positive accuracy.

{\noindent\bf Output:} $d$, a dual lowpass filter with positive accuracy that is biorthogonal to $h$.

{\noindent\bf Step 1:} Find a lowpass filter $f$ that is biorthogonal to $h$.

{\noindent\bf Step 2:} Choose a lowpass filter $g$ with positive accuracy.

{\noindent\bf Step 3:} Compute the polyphase representation ${\tt F}(z)$, ${\tt G}(z)$ (as column vectors), and ${\tt H}(z)$ (as a row vector) of $f$, $g$, and $h$.

{\noindent\bf Step 4:} Set $d$ as the filter whose polyphase representation is ${\tt G}(z)+{\tt F}(z)(1-{\tt H}(z){\tt G}(z))$.

\medskip
We refer to the above as the committee algorithm as the filters $f$ and $g$ (and the input $h$) work together to build the dual of $h$.
The filter $f$ in Step 1 can be found, for example, by using the technique of {G}r{\"{o}}bner bases \cite{P1}.
There are many possible choices for the filter $g$ in Step 2, and one possibility is to take the input $h$ as $g$. Once specific $f$ and $g$ are chosen in Steps 1 and 2, the dual filter $d$ in Step 4 is uniquely determined.

%Although Theorem~1 leads to a wavelet FB associated with the above pair $(h,d)$ of lowpass filters, the wavelet FB associated with the pair $(h,d)$ is not unique in general.

We recall that for a highpass filter $l$, the number of zeros of
$L(e^{i\omega})$, the Fourier transform of $l$, at $\omega=0$ is
referred to as the {\it number of (discrete) vanishing moments} of the
filter $l$ \cite{CHR}.  
It is well-known (see, for example, \cite{CHR}) that for a non-redundant FIR wavelet FB associated with the above pair of lowpass filters $h$ and $d$ with accuracy $\alpha$ and $\beta$ respectively, the synthesis wavelet filters have at least $\alpha$ vanishing moments, and the analysis wavelet filters have at least $\beta$ vanishing moments. Since $d$ is determined by $f$, $g$ and $h$, expressing $\beta$ in terms of the properties of more attainable $f$, $g$ and $h$ is desirable if possible. The next theorem shows that it can be done in a simple way. We omit its proof as it can be obtained by scrutinizing the proof of Theorem 2 in \cite{Hur}. 
 
\begin{theorem}
Let $\beta_1$ and $\beta_2$ be positive integers. 
Suppose that the accuracy of $g$ is $\beta_1$, and 
the number of zeros of $q-H(e^{i\omega})G(e^{i\omega})$ at $\ome=0$ is
$\beta_2$, where $H(e^{i\omega})$ and $G(e^{i\omega})$ are Fourier transforms of $h$ and $g$, respectively.
Then the dual filter $d$ in the committee algorithm has at least 
$\min\{\beta_1,\beta_2\}$ accuracy, hence the analysis wavelet filters in the wavelet FB generated by the pair $(h,d)$ have at least   $\min\{\beta_1,\beta_2\}$ vanishing moments.
\end{theorem}

We finish this section by illustrating our findings using an example.

{\noindent\bf Example: }Let us consider the 1-D dyadic case with the Gaussian filters suggested by Burt and Adelson \cite{BA}. We assume that we are given the lowpass filter $h$ of the form $h=\sqrt{2}\,\tilde{h}$ where
$$\tilde{h}(0)=a,\quad \tilde{h}(-1)=\tilde{h}(1)={1\over 4},\quad \tilde{h}(-2)=\tilde{h}(2)={1\over 4}-{a\over 2},$$
where $a$ is a scalar parameter.
For this lowpass filter, a simple computation shows that, as long as $a\ne 0.25$, the Laurent polynomials 
$f_0(z)=\sqrt{2}\,{1\over 4a-1}$ and $f_1(z)=\sqrt{2}\,{2a-1\over 4a-1}(z+1)$
satisfy the condition (\ref{eq:coprime}). Below we assume that $a\ne 0.25$. Although the filter $f$ associated with the polyphase representation $[f_0(z), f_1(z)]^T$ is always lowpass, its accuracy is zero except the case $a=0.75$, hence it cannot be used as the dual lowpass filter for constructing wavelet FBs directly unless $a=0.75$. However in the committee algorithm, another lowpass filter $g$ that is responsible for the accuracy can be utilized. For simplicity, we choose $g=h$. In this case $\min\{\beta_1,\beta_2\}$ in Theorem 2 is two for every $a\ne 0.25$. For example, for the choice of $a=0.6$, the dual lowpass filter $d$ obtained by the committee algorithm has eleven taps. In this case, both $h$ and $d$ generate the scaling functions that are in $L^2(\RR)$, whose graphs are shown in Figure 1. We note that the choice of $a=0.5$ leads to the interpolatory filter $h$, which is studied in \cite{Hur}.

\begin{figure}
\centering
\includegraphics[width=3.5in]{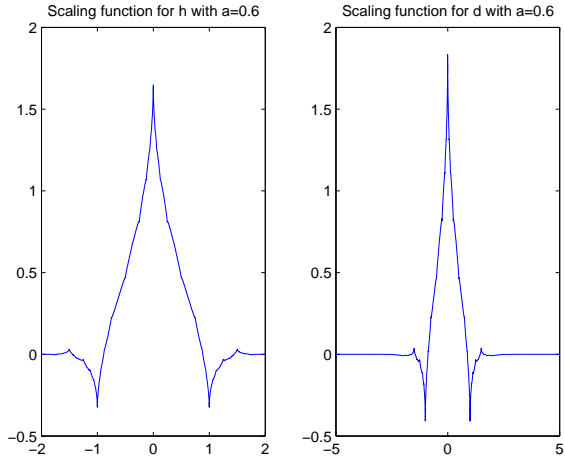}
\caption{Scaling functions associated with the filters $h$ and $d$ for $a=0.6$ in the example.}
\end{figure}
      
\section{Summary}        
In this paper we presented a new technique called the committee algorithm, which suggests an easy way to construct a dual lowpass filter of a given lowpass filter $h$.
Finding a dual lowpass filter $d$ that is biorthogonal to $h$ and has positive accuracy at the same time is a crucial step in constructing non-redundant wavelet FBs, but it is not easy.
In our method, one needs to find two lowpass filters separately (one responsible for biorthogonality and the other for positive accuracy), which is much easier than finding one lowpass filter that satisfies both.
Then the dual lowpass filter $d$ can be obtained in a straightforward manner.
\bibliographystyle{IEEEbib}
\bibliography{IEEEabrv,mybibfile}
\end{document}